\documentclass[journal]{IEEEtran}
\usepackage{graphicx}
\usepackage{amsmath}
\usepackage{amsfonts}
\usepackage{amssymb}
\usepackage{amsthm}
\usepackage{mathtools}
\usepackage{enumerate}
\usepackage{lscape}
\usepackage{longtable}
\usepackage{rotating}
\usepackage{multirow}
\usepackage{url}
\usepackage{array}          % For \hrule and \vrule governing width of table lines (ref: https://tex.stackexchange.com/questions/256731/changing-width-of-table-lines)
\usepackage{caption}   % Apparently this package is needed for the `subfigure` environment. However, it gives warning when being compiled by the IEEE test compilation tool: https://latexqc.ieee.org/report
\usepackage{subcaption}
\usepackage{rotating}
\usepackage{comment}
\usepackage[utf8]{inputenc}
\usepackage{amssymb}
\usepackage{amsfonts}
\usepackage{tikz}
\usepackage{natbib} 

\usepackage{color}
\usepackage{hyperref}
\usepackage[shortlabels]{enumitem}
\usepackage{ulem}       % For \sout (strike out)

\def\A{\mathcal{A}}
\def\C{\mathcal{C}}
\def\E{\mathbb{E}}

\def\P{\mathbb{P}}
\def\ind{\mathbb{I}}
\def\FV{FV }%{Fleming-Viot particles}

\definecolor{darkgreen}{rgb}{0.0, 0.5, 0.0}

\title{Multistage Stochastic Programming for Rare Event Risk Mitigation in Power Systems Management}
\author{Daniel Mastropietro, Vyacheslav Kungurtsev
\thanks{D. Mastropietro and V. Kungurtsev are with the Department of Computer Science at the Czech Technical University, Praha, Czech Republic}
}
\begin{document}

\maketitle

\footnote{This work was submitted in April 2026 to the IEEE for possible publication. Copyright may be transferred without notice, after which this version may no longer be accessible.}

\begin{abstract}
High intermittent renewable penetration in the energy mix presents challenges in robustness for the management of power systems' operation. If a tail realization of the distribution of weather yields a prolonged period of time during which solar irradiation and wind speed are insufficient for satisfying energy demand, then it becomes critical to ramp up the generation of conventional power plants with adequate foresight. This event trigger is costly, and inaccurate forecasting can either be wasteful or yield catastrophic undersupply. This encourages particular attention to accurate modeling of the noise and the resulting dynamics within the aforementioned scenario. In this work we present a method for rare event-aware control of power systems
%that uses a Fleming-Viot particle approach to manage the control forecasts of such tail events to manage the control forecasts of such tail events.
using multi-stage scenario-based stochastic programming.
A Fleming-Viot particle approach is used to bias the scenario generation towards rare realizations of very low wind power, in order to obtain a cost-effective control of conventional power plants that is robust under prolonged renewable energy shortfalls.
\end{abstract}

\begin{IEEEkeywords}
% Source: https://ieee-org.widen.net/s/dr9nl9n8kw/ieee-taxonomy_jan_2026_pdf
Fleming-Viot, renewable energy, robustness, wind energy
\end{IEEEkeywords}

\section{Introduction}
Management of power systems involves timing the generation of energy across various producers in the energy mix in a manner that satisfies the energy demand while at the same time seeking to minimize cost and maximize generation from renewables. Renewables such as wind and solar are subject to intermittency - their generating capacity on any day can vary depending on solar irradiation and wind speed. At the same time, conventional plants require time to start up before they can generate power. Thus, as renewables' prevalence in the energy mix increases, the ability to accurately predict when to fire up conventional generation becomes increasingly critical to ensure the stable and adequate provision of energy.

In this paper we model low renewable generation due to weather as a problem in rare event simulation. We extend recent work~\cite{tapia2025electricitymarketclearingextremeevents} that applies an approach that incorporates the large deviations principle towards a multistage finite time horizon framework. We use the Fleming-Viot particle scheme (FV) \cite{AFG} to preferentially forecast scenarios that are of interest for robust operation. The operating model uses market clearing in energy with wind and conventional coal plants. We develop a scheme that develops accurate forecasts of low wind generation scenarios and starts the operation of conventional generation as needed to satisfy the projected shortfall in demand satisfaction. Both conventional plants require multiple time stages to activate, and so it becomes necessary to predict in advance whether their generation is needed to forestall shortfall. 

The application of FV and use of rare event and risk optimization metrics are associated with the consideration of the critical contemporary problem of renewable integration. Because solar and wind are intermittent, it is important to have
%alternatively
alternatives available to provide energy supply in the case of low wind speed and solar irradiation, with conventional coal plants the most reasonable choice with current technology. In particular, events described as dunkelflaute, the German phrase meaning ``dark lull'', i.e., darkness and still air, presents a serious risk for energy systems with a high penetration of renewable sources (see, e.g.~\cite{kittel2024coping}). If these events are not sufficiently forecasted and mitigated for in advance, a real hazard of comprehensive grid failure arises. Thus,  while we include the costs of ramp up and shut down of the conventional plants in the optimization, the priority as far as the forecasting ---that is the scenarios generated--- and the optimization criterion are concerned, is towards rare events of low and exceptionally low wind generation.

Below we review the (brief) set of works considering similar approaches. Afterwards, in Section~\ref{sec:problem} we present the multistage stochastic programming problem and the time series model for the noise. Section~\ref{sec:methods} describes our method and implementation details, including a description of the FV algorithm for rare event sampling. Section~\ref{sec:results} presents our numerical results. We end with some concluding comments in Section~\ref{sec:conc}

\subsection{Related work}
The necessity of energy supply security in the presence of the risk of dunkelflaute periods is discussed in detail in~\cite{pelka2018impact}. A robust worst-case energy management strategy is presented in~\cite{bernecker2025adaptive}. This approach is known to be excessively conservative, pushing up average costs. Battery storage is one mitigation strategy that has been investigated~\cite{lee2025interval}, however the technology for sufficient storage capacity during prolonged low renewable generation is not yet available in practice. The well known and widespread Stochastic Dual Dynamic Programming (SDDP) approach for adaptively ramping up conventional coal-fired plants in the presence of insufficient renewable generation is presented in~\cite{li2023retrofit}. In this paper, we consider a tight chance constraint preventing undersupply, which presents a nonconvex problem that is unsuitable for SDDP. Another strategy is to use Distributionally Robust optimization, which is incorporated in a hybrid long term planning and real-time dispatch stochastic programming operation in~\cite{huang2021two}.

This work builds on a recent paper that incorporates probabilistic tail modeling explicitly. The work in \cite{tapia2025electricitymarketclearingextremeevents} presents electricity market clearing in the presence of extreme events, using large deviation theory (LDT) to estimate the probability of those deviations.

The authors present two chance-constrained optimization problems that aim at minimizing the expected power generation cost, where the power generation consists of a scheduled generation plus a random component consisting of a reserved flexible capacity against the occurrence of rare events that put at risk power demand satisfaction.
The two models use LDT to propose a control policy for distributing the reserve capacity into the different power generators that can controllably inject power into the network. The models are distinguished by the use of weights that differentiate two types of risk, one far away from the rare event and one closer to it.
In this work we focus on one type of risk but extend their chance-constrained formulation to a \textit{multistage} stochastic optimization problem.

\section{Problem description}
\label{sec:problem}
The specific problem considered for illustration consists of a power system supplied by two types of energy sources: renewable (wind farms) and coal.
The renewable source is considered a stochastic process and coal is a backup source controllable by a system operator.
The general policy of the system operator is the following: renewable sources are used as long as they can satisfy the demand; when this is not possible, coal plants are expected to provide the remaining energy requirement. However, this needs to be forecasted and decided preemptively, since there is a period of time that must elapse for a coal plant to turn online or shut off.
%Control decisions are only possible on the coal plants, which
This is modeled by the coal plants operating, at any given time, in one of the following four possible states: idle, starting, operating, stopping.

The goal of the optimization problem is to design a multi-stage decision plan for the state of the coal plants. For simplicity, we consider that there is just one coal plant that: (i) takes one unit of time (stage) to start and one unit of time to stop, and (ii) has sufficient capacity to provide power for any realization in the support of the energy demand distribution.

\bigskip
\noindent \textbf{Model}

%\todo{Use $x$ for states, which depend on $\xi$ (random) (Shapiro); use $z$ or $u$ for the decisions. $\Omega$ is the space of possible outcomes in a finite horizon of length $H$.}

We follow \cite{tapia2025electricitymarketclearingextremeevents} to write the problem as a chance-constrained multi-stage optimization problem, where a planning horizon of length $H$ defines $H$ decision stages. For each stage $h, 0 \leq h \leq H-1$, we let $D_h$ and $W_h$ be two non-negative random variables representing, respectively, the energy demand and the generated wind power.
%$S$ be the renewable source with large availability uncertainty, in our case, the swarm of wind farms. Given the known wind power provided at the current stage $0$ and a planning horizon $H > 0$, we let
%$W_h$ the generated wind power,
%by the swarm at stage $h$,
We also let $\kappa_h \in \{1, 2, 3, 4\}$ the actionable state of the coal plant at the respective stage $h$, where $1$ means ``idle'', $2$ means ``starting'', $3$ means ``operating'', and $4$ means ``stopping''. We finally define the decision variable to optimize at each stage $h$ as a four-dimensional vector, $\mathbf{x}_h := (x^j_h)_{j=1,2,3,4}$, with coordinates given by $0/1$ variables indicating each of the four possible coal plant states, namely $x^j_h := \ind_{\kappa_h=j}$ for $0 \leq h \leq H-1$, where $\ind$ is the indicator function. We also define $\mathbf{x} := ((x^j_h)_{j=1,2,3,4})_{h=0, \dots, H-1} \in \{0, 1\}^{4H}$ as the vector of all binary decision variables across the planning horizon of length $H$. Linear constraints enforce physical feasibility of stage-to-stage decisions, e.g., a plant cannot be off at some stage $h$ and operational at stage $h+1$.

%The problem considered is the multistage chance-constrained economic dispatch with large deviation theory \todo{(LDT-CC-ED) \cite{tapia2025electricitymarketclearingextremeevents}}, where
The objective is to minimize the future expected total power generation cost from non-renewable sources, while limiting the rate of violations to the probability (or chance) constraint that the demand is satisfied.
%due to the uncertainty renewable power (the wind farms in this case).
%in the planning horizon of length $H$.
%By approximating the future expected total power as the expected total power over the planning horizon $H$ and by letting $\Omega$ be the space of all possible wind power outcomes in said horizon, the optimization problem can be written as:
We let $\boldsymbol\Omega := (\Omega, \mathcal{F},\mathbb{P})$ be the probability measure space for the demand and wind power stochastic processes, with $\Omega$ the space of all possible outcomes of demand and wind power across all $H$ stages, $\mathcal{F} := \{\mathcal{F}_h\}_{h=0, \dots, H-1}$ the collection of event filtrations across stages, and $\mathbb{P}$ the probability measure of events defined as subsets of $\Omega$. We also write $\Omega_h$ the space of possible outcomes at stage $h$, and $\omega_h$ a particular outcome in $\Omega_h$. The decision $\mathbf{x}_h$ is a function of the outcomes \textit{up to stage} $h$, i.e. $\mathbf{x}_h := \mathbf{x}_h(\omega_{[h]})$, where $\omega_{[h]} := (\omega_0, \omega_1, \dots, \omega_{H-1})$, which is also referred to as the \textit{non-anticipativity constraint}. We also let $f_h$ the cost associated to the coal plant at stage $h$ ---which is a function of the demand, the wind power, and the coal plant state at stage $h$, as described next.
Given a small probability $\epsilon > 0$, we write the chance-constrained optimization problem as:

\begin{subequations}
\label{equ:optim_problem}
\begin{alignat}{2} % This implies a 4-column layout with rlrl alignment, so we leave the first right alighment empty then add a & to signify the first left alignment and leave the second right alignment again meaning that the && marks the second left aligned column. Ref: https://tex.stackexchange.com/questions/200502/align-environment-align-on-the-left-side for LEFT ALIGNMENT
    \label{equ:optim_problem_objective}
    & \min_{\mathbf{x} \in \{0,1\}^{4H}} \E_{\boldsymbol{\Omega}} \left[\sum_{h=0}^{H-1} f\big(\mathbf{x}_h(\omega_{[h]}), D_h(\omega_h), W_h(\omega_h)\big) \right] \\
    & \nonumber \text{s.t.} \\
    \label{equ:optim_problem_decision_variables}
    & \mathbf{x} := \{\mathbf{x}_h(\omega_{[h]})\}_{h=0, \dots, H-1} \\
    & \nonumber \P_{\boldsymbol{\Omega}}\big[p_h x^3_h \geq D_h - W_h, \\
    \label{equ:optim_problem_constraint_power}
    & \qquad \qquad \qquad 0 \leq h \leq H-1 \big] \geq 1 - \epsilon \quad \text{\footnotesize (chance constraint)} \\
    & \nonumber \text{for } 0 \leq h \leq H-1: \\
    & \text{ }x^1_h, x^2_h, x^3_h, x^4_h \in \{0, 1\} \quad \quad \text{\footnotesize (binary constraint)} \\
    & \text{ }x^1_h + x^2_h + x^3_h + x^4_h = 1 \quad \text{\footnotesize (single-state constraint)} \\
    & \nonumber \text{for } 1 \leq h \leq H-1: \\
    & \nonumber \text{\footnotesize (temporal constraints on plant operation at consecutive stages)} \\
    & \text{ }x^1_{h-1} \leq x^1_h + x^2_h \text{  
 \footnotesize (IDLE } => \text{\footnotesize IDLE or START.)} \\
    & \text{ }x^2_{h-1} \leq x^3_h + x^4_h \text{  
 \footnotesize (START. } => \text{\footnotesize OPER. or STOP.)} \\
    & \text{ }x^3_{h-1} \leq x^3_h + x^4_h \text{  
 \footnotesize (OPER. } => \text{\footnotesize OPER. or STOP.)} \\
    & \text{ }x^4_{h-1} \leq x^1_h + x^2_h \text{  
 \footnotesize (STOP. } => \text{\footnotesize IDLE or START.)},
\end{alignat}
\end{subequations}
where $p_h := \min\{p_{\max}, [D_h - W_h]_+\}$ is the power provided by the coal plant of capacity $p_{\max}$, if operating at stage $h$.
%$p_h = \min([D_h - \hat{W}_h]_+, p_{\max})$ and
%$p_{\max}$ is its power capacity.
%\todo{Note that, if $p_{\max}$ is smaller than the power shortfall, $[D_h - W_h]_+$, constraint \eqref{equ:optim_problem_constraint_power} will not be satisfied and the problem will be infeasible.}
Given, respectively, $b_2$, $b_3$, $b_4$ the starting, operating, and stopping costs of the plant, and $b$ the cost per unit of generated power, the cost $f_h$ associated to the coal plant at stage $h$ is computed as:
\begin{equation}
\label{equ:cost}
    f(\mathbf{x}_h(\omega_{[h]}), D_h(\omega_h), W_h(\omega_h)) := b_2 x^2_h + (b_3 + b p_h) x^3_h + b_4 x^4_h.
\end{equation}
%with $p_h = \min\{p_{\max}, [D_h - W_h]_+\}$, when $x^3_h = 1$, i.e. when the plant is ``operating''.

%We note that the setup of the problem assumes that the coal plant may be already started or operating before the very first stage $h = 0$, therefore no constraint is imposed on the coal plant being either idle or starting at the initial stage. \textcolor{red}{I thought you always start stage 0 to have enough wind generation to cover demand and the coal plant off?} \dan{We discussed the two approaches and I thought you preferred the approach that assumes an ``ongoing'' operation of the plant. If we want to assume the idle or start state at stage $0$, we should add a constraint stating this, i.e. $x^1_0 + x^2_0 = 1$ (which is currently commented out). However, I am not using this constraint in my implementation.}

\bigskip

\section{Methodology}
\label{sec:methods}
%\todo{Write the problem when we discretize the state space and define the scenarios which correspond to a finite set of samples. use Shapiro's book as reference, where they define e.g. a $\xi$ variable. \newline
%Change the chance constraints to almost surely. See also when Shapiro talk about chance constraint and almost sure constraints. \newline
%}

%\todo{We should say here that we make the chance constraint be satisfied for all the scenarios, making the solution more conservative than the original problem. }

We start by noting that problem \eqref{equ:optim_problem} is an infinite dimensional optimization problem because, as per \eqref{equ:optim_problem_decision_variables}, decision variables to optimize, $\mathbf{x}_h$, are \textit{functions} of real-valued stochastic processes $D_t$ and $W_t$ for $t \leq h$, and are therefore measure-valued. This presents an optimization problem of infinite dimension unless $D_t$ and $W_t$ have a finite number of possible realizations.
%, which take values on an infinite dimensional space.

Therefore, following \cite[Chapter 3]{shapiro2021lectures}, we propose to approximate the problem numerically using scenarios. This technique allows computing a sample average approximation of the objective function to obtain an approximate solution of the original infinite dimensional optimization problem.
The quality of the approximation and its convergence properties as the number of scenarios goes to infinity are described in \cite[Section 5.8]{shapiro2021lectures}.
%which discusses the statistical properties of the sample average approximation estimators of the objective expectation.

In the scenario approach, we add another layer of approximation by replacing all chance constraints in the original problem with hard constraints, making the approximated solution more conservative than without these replacements. Standard concentration inequalities (e.g.~\cite{boucheron2003concentration}) provide justification that with a sufficient number of scenarios, a constraint holding for all scenarios implies that the chance constraint holds. 
%that should be satisfied by every considered scenario. This makes the approximated solution more conservative than otherwise keeping the chance constraints.

Thus, we solve the optimization problem obtained from \eqref{equ:optim_problem} after performing the following substitutions in the given two expressions:
\begin{align}
    \label{equ:optim_problem_objective_estimated} \nonumber \text{\eqref{equ:optim_problem_objective}} & \rightarrow \min_{\mathbf{x} \in \{0,1\}^{4H}} \quad \hat{\E}_{\boldsymbol{\Omega}} \left[\sum_{h=0}^{H-1} f_h (\omega_{[h]}) \right] := \\
    & \qquad \qquad \qquad \min_{\mathbf{x} \in \{0,1\}^{4V}} \frac{1}{S} \sum_{i=1}^{S} \sum_{h=0}^{H-1} f_h (\omega_{i,[h]}) \\
    \label{equ:chance_constraint_as_hard_constraint} \nonumber \text{\eqref{equ:optim_problem_constraint_power}} & \rightarrow p_h(\omega_{i,h})x^3_h(\omega_{i,[h]}) = p_h(\omega_{i,h}) \quad \forall i \in [S] \\
    & \quad \text{with } p_h(\omega_{i,h}) := \min\{p_{\max}, [D_h(\omega_{i,h}) - W_h(\omega_{i,h})]_+\},
\end{align}
where $S$ is the number of generated scenarios, $V$ is the number of nodes in the scenario tree, $\omega_{i,h}$ is the outcome of the stochastic processes at stage $h$ in scenario $i$, and $\omega_{i,[h]}$ is its history up to stage $h$.
Note that the new form of constraint \eqref{equ:optim_problem_constraint_power} enforces the operation of the plant ($x^3_h = 1$) when the system requires a positive power generation $p_h$.

%The optimization problem in \eqref{equ:optim_problem} has the form of the linear multi-stage optimization problem presented in \cite[Chapter 3]{shapiro2021lectures}. Due to the infinite state space on which the stochastic processes $D_h$ and $W_h$ are defined, the optimization problem has infinite dimension. Its solution, however, can be approximated numerically using scenario trees, as described in the aforementioned chapter.

In the scenario approach, the stochastic processes involved in the problem are considered to take a finite number of values at each stage by way of a branching procedure that generates a scenario tree. Taking the $(W_h)_{0 \leq h \leq H-1}$ stochastic process for illustration and having defined a suitable data generation process, we construct a scenario tree of order $B$ as follows: we place the known $W_0$ value at the root of the tree and generate $B$ finitely many possible values for $W_1$ to create $B$ nodes at stage $1$; for each of these nodes, we generate $B$ possible values for $W_2$, independently of $W_1$, thus creating $B^2$ nodes at stage 2; and so forth until reaching stage $H-1$.
This gives rise to a tree of order $B$ with $1 + B + B^2 + ... + B^{H-1} = \frac{1 - B^H}{1 - B}$ nodes, and $S := B^{H-1}$ scenarios (given by the number of leaf nodes).

The linear optimization problem is then written as a ``giant'' linear programming problem by associating a different optimization variable to each node in the tree.
%The problem is completed by adding the non-anticipativity constraints, that ensure that the decision variable at stage $h$, $\mathbf{x}_h$.
With this formulation, the decision variable at each node $k$ of stage $h$, $\mathbf{x}_{h,k}$, depends on the realization of the stochastic processes up to stage $h$ leading to node $k$, namely on the particular realizations of $D_{[h]} := (D_0, D_1, \dots, D_h)$ and $W_{[h]} := (W_0, W_1, \dots, W_h)$.

In our implementation, we assume the demand $D_h$ has low variability and the wind power $W_h$ is a process with high variability, whose impact on the system we would like to mitigate by appropriately controlling the coal plant operation.
Thus, for illustration of the methodology, we model $D_h$ as a constant value and focus the scenario tree generation on the $W_h$ process.

The approach for obtaining robust control strategies of the coal plant, capable of compensating wind power supply shortfalls, is to \textit{bias} the data generation process towards rarely observed realizations. We achieve this by (i) modeling the \textit{change} in wind power between consecutive stages, and (ii) splitting the $B$ branches at each node into two separate data generation processes: one that models ``normal'' changes and one that models ``rare'' changes. 
%alongside scenarios with typical wind power evolutions, a representative number of scenarios with low-valued rare deviations.
%in order to also control for rare but impactful wind power shortfall situations.
%We achieve this by modeling the \textit{change} in wind power between consecutive stages, and oversampling negative change realizations.
%Additionally, both data generation processes model the \textit{change} in wind power between consecutive stages.
The normal changes are generated using an autoregressive (AR) model, while rare negative changes are generated using Fleming-Viot particle system techniques adapted from \cite{mastropietro2025questa}, as described next.

%in combination with weighted-risk strategies that prioritize scenarios with higher risks \todo{(e.g. coming from a Fleming-Viot discovery of these risks. This means: simulate the scenario-generation process with Fleming-Viot in order to ACTUALLY GENERATE scenarios that imply high risks.)}.

%\subsection{Scenario generation}
%\todo{TBC. Probably section not needed.}
\bigskip
\noindent\textbf{Fleming-Viot particle systems}
\label{sec:methods_fv}
%\todo{Summarize and refer to the literature for more details.}

%\todo{It allows having a good estimate of the expectation when we bias the samples in the scenario tree.}

We use Fleming-Viot particle systems (FV) as a technique to estimate the probability of unlikely changes in wind power supply. Using a regenerative argument, the technique is suitable for providing a reliable estimate of rare event probabilities as presented in \cite{mastropietro2025questa} to the context of queueing systems. We now provide a summary of the estimation methodology.

Given a continuous-time ergodic Markov process $\{\zeta_t\}_{t \geq 0}$, the method leverages the regenerative characterization of the stationary probability of a set $\C$ given by \cite[Ch.~6]{asmussen}:
\begin{equation}
\label{equ:proba_stationary}
    p(\C) = \frac{\E_A \int_{0}^{T_A} \ind_{\{\zeta_t \in \C\}} dt}{\E_A T_A},
\end{equation}
where $A$ is an event defining a regenerative cycle and $T_A$ is the regenerative cycle time.
%\FV particle systems were introduced by \cite{AFG} to consistently estimate quasi-stationary distributions of Markov processes with absorption.
The \FV estimation approach builds on the now matured use of Fleming-Viot particle systems to consistently estimate quasi-stationary distributions of Markov processes with absorption \cite{AFG}.
As shown in \cite{mastropietro2025questa}, the numerator in expression~\eqref{equ:proba_stationary} can be written in terms of an approximation of the quasi-stationary distribution, $\nu$, of a process that is absorbed when visiting a set $\A$, directly linked to the event $A$ defined above.
%defining the regenerative cycle of the aforementioned Markov process $\zeta_t$.
Thus, event $A$ in \eqref{equ:proba_stationary} becomes the \textit{absorption event}, and $\A$ is called the \textit{absorption set}.

The \FV system consists of $N$ copies (called particles) of the $\zeta_t$ process constrained to evolve outside $\A$ by an appropriate restarting mechanism of a particle touching it. By defining the absorption set $\A$ as a set of frequently visited states by $\zeta_t$, the \FV system is able to collect further samples from regions rarely visited by $\zeta_t$.
Using the convergence guarantees of the FV estimator of the quasi-stationary distribution $\nu$ as the number of particles $N \to +\infty$ \cite{AFG}, a consistent estimator of the stationary probability $p(\C)$ in \eqref{equ:proba_stationary} can be constructed, as long as $\C$ is \textit{outside} $\A$ \cite{mastropietro2025questa}.
%(see \cite{mastropietro2025questa} for details), 
%enabling a more accurate probability estimation than a Monte Carlo exploration by $\M$.
When the set of interest $\C$ is a low-probability set, this FV estimator is expected to provide a more accurate estimation of $p(\C)$ than using \eqref{equ:proba_stationary} on samples drawn from a conventional Monte Carlo discretization of $\zeta_t$.

\bigskip
\noindent \textbf{Adaptation to scenario tree generation}
\medskip

We propose using \FV to generate rare negative changes in the wind power process that are large enough to put power demand satisfaction at risk. We base the rare generation process on reliable estimations of the probability of selected large negative changes in the wind power process.

\begin{figure}[!ht]
    \centering
    \begin{subfigure}[b]{0.45\textwidth}
        \includegraphics[width=0.9\linewidth]{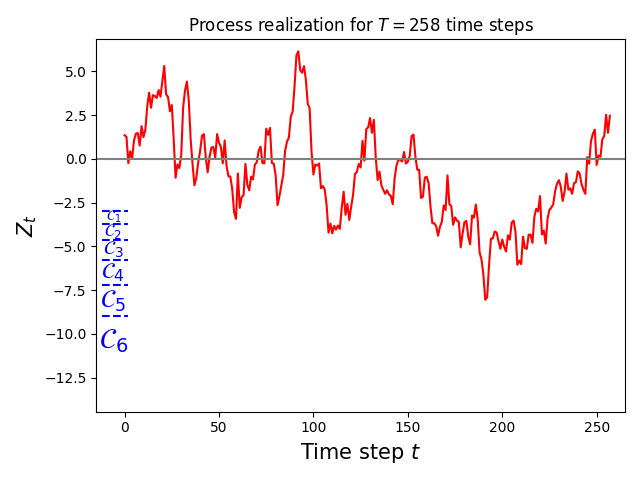}
        \caption{Process sample realization}
    \end{subfigure}
    \begin{subfigure}[b]{0.45\textwidth}
        \includegraphics[width=0.9\linewidth]{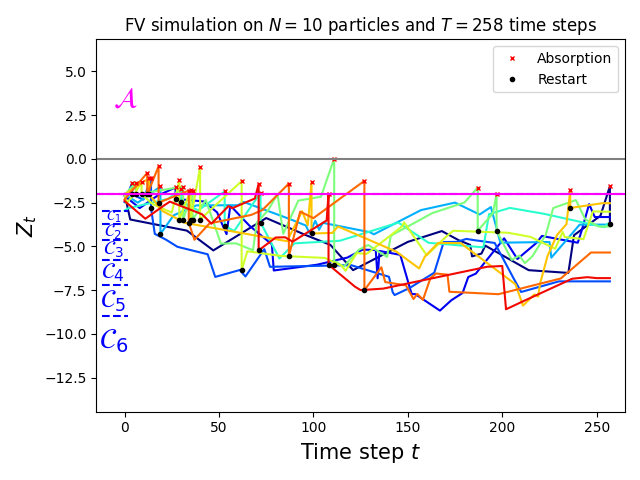}
        \caption{Fleming-Viot (FV) simulation with $10$ particles}
    \end{subfigure}    
    \caption{(a) Typical realization of the $Z_t$ process defined as a positively-correlated AR(2) process, $Z_t = 0.90 Z_{t-1} + 0.05 Z_{t-2} + \text{error term}$ (see Section~\ref{sec:scenario_generation}).
    (b) Example of an FV simulation with $N = 10$ particles, where each color represents a different particle. The absorption set $\A$ (lower bounded by the horizontal magenta line) and the pre-defined disjoint intervals $C_k$ (bounded by the horizontal blue lines) are overlaid.
    Note that the time steps at which each particle is updated vary across particles, since the discrete-time FV dynamics chooses uniformly at random the next particle to update at each time step.
    \\
    The two plots share the same horizontal and vertical axes as well as the same number of simulation steps $T$, so that they can be compared fairly in terms of exploration budget.}
\label{fig:fv_and_mc}
\end{figure}

To this end, we define process $\{Z_t\}_{t \in \mathbb{N}}$ as a time series representing the change in wind power between two consecutive stages,
%in the planning horizon,
that is $Z_t := W_t - W_{t-1}, t \geq 1$.
Since large \textit{negative} changes are of interest,
%that may hinder the satisfaction of power demand,
we define the absorption set $\A$ of the \FV system as the semi-open interval $[-a, +\infty)$ for some positive $a$ for which the probability that $Z_t$ is in $\A$ is considered large (which makes $\A$ uninteresting for further exploration).
In order to define a formal Markov process for performing Fleming-Viot, we augment the state to include the appropriate recent history of the time series model for $Z_t$\footnote{If $Z_t$ is an ARMA$(p, q)$ process, the Markov state is defined of dimension $p$ as $\mathbf{U}_t := (Z_{t-p+1}, \dots, Z_{t})$, as this allows the next state value to be a function of the current state value. Recall that the MA process is given by a linear combination of an uncorrelated Gaussian process $\epsilon_t$, whose values can be generated independently of the $Z_t$ values and of past $\epsilon_t$ values.}.
We then use the extension of \FV to \textit{discrete-time} stochastic processes \cite{mastropietro2024discretetimefv} (where just one particle, selected uniformly at random among all particles, is updated at every time step) to estimate the probability of large negative values of $Z_t$
%(i.e. large negative changes in wind power $W_t$),
belonging to predefined disjoint intervals $\{\C_k\}_{k = 1, \dots, K}$ into which the interval $\C := (-\infty, -c)$, for some $c > a$, is partitioned\footnote{We are abusing notation when the model for $Z_t$ is an ARMA$(p, q)$ process with $p > 1$ because the $\C$ interval is \textit{not} the same as the $\C$ set in \eqref{equ:proba_stationary}, as the latter is a multidimensional state. The $\C$ set in that expression is actually the $\C$ interval union $\mathbb{R}^{p-1}$, and the FV estimated probability is the probability of the set, not the interval, which however is considered an approximation of the probability of the interval for the data generation process.}. The value of $c$ is chosen so that the $\C$ interval is considered a low probability set (see Figure~\ref{fig:fv_and_mc}(b) for an example). Finally, the \FV system is used to estimate the probabilities of the different $\C_k$ disjoint intervals, where the key input is the proportion of particles in each interval $\C_k$ at each time step $t$ for the duration of the FV simulation (see \cite{mastropietro2025questa,mastropietro2024discretetimefv} for details). These estimated probabilities are used in the scenario generation process for obtaining sensible rare negative changes in wind power.
The next section describes the details of the scenario generation process, including sensible choices of $a$ and $c$ given the time series model proposed for $Z_t$.
%Figure~\ref{fig:fv_and_mc}(b) shows a typical excursion of an FV particle system in that context.
%is illustrated in Figure~\ref{fig:fv_and_mc}(b) for the $Z_t$ process defined in Section~\ref{sec:scenario_generation}. The choices of $a$ and $c$ as well as the scenario generation details are described in that same section.

Figure~\ref{fig:fv_and_mc} illustrates the effect of the FV particle system on biasing exploration towards sets rarely visited by the original Markov process: subfigure (a) depicts a typical realization of the wind power change process $Z_t$, and subfigure(b) shows an FV excursion with a system of $10$ particles, with the choices of $\A$ and $\C_k$ defined in Section~\ref{sec:scenario_generation}.
Thanks to the absorption-and-restart mechanism of the FV particle system, this trajectory spends more time in areas around the $\C_k$ intervals than the unconstrained realization of $Z_t$ depicted in subfigure (a), which, using the consistency of the FV estimator, enables a more reliable estimation than plain Monte Carlo of the probability of those rarely visited sets.
%an ensemble of intervals having lower steady-state probability than the frequently visited absorption set $\A$.

\section{Experimental setup and results}
\label{sec:results}
We evaluate our method on an optimization problem with $H=5$ stages of the form given in \eqref{equ:optim_problem} ---with \eqref{equ:optim_problem_objective} and \eqref{equ:optim_problem_constraint_power} replaced by \eqref{equ:optim_problem_objective_estimated} and \eqref{equ:chance_constraint_as_hard_constraint}--- using the scenario tree approach described in Section~\ref{sec:methods}, with $B = 20$ branches per node. This gives rise to a tree with $168\,421$ nodes, providing $S := 20^{H-1} = 160\,000$ scenarios. %with a conservative value for \todo{$\epsilon = ...$ (TODO, based on large deviations analysis --> summarize this with a small paragraph of what should be done --conditioned on the plant being off, probability that the wind power goes below the demand.)}.

The constant power demand at all stages is set to $D_h = 6$ GW ($h = 0, \dots, 4)$ and the initial wind power supply is set at $W_0 = 10$ GW, so that demand is initially satisfied by renewable sources and the coal plant can be assumed to be in ``idle'' or ``starting'' state.
The capacity of the coal plant is set at a large value of $p_{\max} = 400$ GW to ensure, given the scenario generation model described next, that the chance constraint is satisfied in every generated scenario.
This condition gives rise to a conservative deterministic version of the chance-constrained optimization problem in \eqref{equ:optim_problem}.

The cost parameters defining the operational cost $f$ of the coal plant in \eqref{equ:cost} are set (in some appropriate cost units) to $b_2 = 3, b_3 = 5, b_4 = 2, b = 20/\text{GW}$.

\subsection{Scenario generation}
\label{sec:scenario_generation}
%\todo{TODO: Consider two different models for $Z_h$... one with positive correlation and one with negative correlation. Goal: given that we are not basing the simulations on an actual model of the wind power supply, show that our approach is more robust than the benchmark in both these opposite scenarios. Note that in each case, the definition of the absorption set $\A$ and the set of interest $\C$ should change.}

%\todo{TODO: Cite the paper about the modeling of wind power which would probably suggest a positive correlated model. --> I checked the source (\href{https://khazna.ku.ac.ae/ws/portalfiles/portal/19110843/Modeling-And-Forecasting-Of-Wind-Speed-Time-Series-In-The-United-Arab-Emirates.pdf}{https://khazna.ku.ac.ae/ws/portalfiles/portal/19110843/Modeling-And-Forecasting-Of-Wind-Speed-Time-Series-In-The-United-Arab-Emirates.pdf}), a Master's thesis that models the wind speed in the UAE. It does NOT give details of the fitted PARIMA model (periodic-ARIMA), i.e. it does not give the actual coefficients of the fitted series and it only shows the ACF of the series, NO PACF is shown in any part of the thesis, where most of the analysis is done on the residuals of the fit but not on the actual series.}

The wind power provided by a wind farm is a nonlinear function of the wind speed \cite{windpower2020review}. However, to simplify implementation and focus the analysis on our proposed optimization approach, we assume a linear relationship between wind power and wind speed. Thus, we generate the values of $W_h$ using this simplified model of the relationship of wind speed and power generation, with the former a stochastic process of the form given in the extensive literature on the topic.
%and model the wind speed, of which there is extensive literature.
%That is, we do not go into the details of taking into account their highly nonlinear relationship, such as the fact that wind mills are stopped if the wind speed goes above a given threshold established for the safe operation of the mill (see \cite{windpower2020review} for details).
%Thus, the wind speed process $W_h$ and the wind speed change process $Z_h$ are also regarded as wind power processes and expressed in gigawatts (GW).

Several works modeling and forecasting wind speed and energy in different regions of the world (e.g. United States, Poland, Baltic Sea, India, United Arab Emirates) propose seasonal ARIMA models with a difference order equal to $d = 1$ (see for instance \cite{windpower2009farima,windpower2020review}).
In our scenario generation, we propose using an ARIMA model for the wind speed, which relates directly to the wind power $W_h$, assumed to represent average hourly power measurements. We model the series as a differenced series of order $d = 1$ and, for simplicity, we don't use any moving average (MA) term nor seasonality components. The latter can be justified by the fact that we will look at very short-term horizons in our analyses (e.g. 5 hours) where we do not expect any observable seasonality effect.
Following the literature, we further assume a positively autocorrelated time series and propose a zero-mean stationary AR(2) model for the differenced wind power series, $Z_h$, with parameters $\phi_1 = 0.90$, $\phi_2 = 0.05$, that is, $Z_h = 0.90 Z_{h-1} + 0.05 Z_{h-2} + \epsilon_h$, where $\epsilon_h$ is a zero-mean Gaussian random process with unit variance, independent of $Z_h$.
%Thus, we model the differenced wind power series, $Z_h$, as an autoregressive (AR) model, which is used as instances of the ``normal'' behavior of wind power change.
%(guaranteed by the unit roots of the characteristic polynomial being smaller than $1$ in absolute value \todo{REF}).
%Realizations of $Z_h$ are
%This model is used as the generative model of the ``normal'' behavior of wind power change.e

\medskip

This generative process is used in two scenario generation contexts, to be evaluated and compared:
\begin{enumerate}[label=\alph*)]
    \item \label{itm:bm} \textbf{Benchmark:} A reference approach where the data generation is governed solely by the $Z_h$ process.

    \item \label{itm:ours} \textbf{Biased:} Our proposed approach where the data generation follows a mixture between the $Z_h$ process and rare changes in wind power, whose magnitude and probability is governed by the \FV system described in Section~\ref{sec:methods}.
    %is biased towards rarely observed realizations, as described in Section~\ref{sec:methods}.
\end{enumerate}
In (\ref{itm:bm}, the realizations of the $Z_h$ process are assigned to all $B$ branches. In (\ref{itm:ours}, the realizations of $Z_h$ are assigned to half of the $B$ branches created at each stage of the scenario tree and the other half is assigned to ``rare'' negative realizations of the wind power change, as described next.

The ``rare'' realizations of wind power supply in (\ref{itm:ours} are based on the probabilities of the $\C_k$ intervals estimated by FV described in Section~\ref{sec:methods_fv},
%are used to generate large negative change realizations in wind power, and assigned to the other half of the $B$ branches at each stage.
%Following are the details for the Fleming-Viot particle system setup and for the data generation process.
and are obtained with the following steps (see left plot in Figure~\ref{fig:sets_and_realizations} for the visualization of these sets and intervals):
\begin{enumerate}
    \item \textbf{Definition of absorption set $\A := [-a, +\infty)$:} Based on the trade-off considerations outlined in \cite{mastropietro2025questa}, the value of $a$, defining the absorption set $\A$ containing frequently visited states, is chosen as $2.0$ standard deviations of the $Z_h$ process, namely as $a = 2.0$. This makes the absorption set of FV be selected as $\A := [-2.0, +\infty)$.
    The absorption event is defined by the current $Z_h$ value touching $\A$ while the previous $Z_{h-1}$ was outside $\A$.
    
    \item \textbf{Definition of the interval of interest $\C := (-\infty, -c)$ and its sub-intervals:} The value of $c$ is chosen as $3.0$ standard deviations, namely as $c = 3.0$, making $\C := (-\infty, -3.0)$.
    This interval is partitioned into $6$ sub-intervals, as follows: $\C$ is first partitioned into two intervals $I_1 := [-9.0, -3.0)$ and $I_2 := (-\infty, -9.0)$; then, $I_1$ is partitioned into $5$ sub-intervals of equal length $\Delta$ in the log-10 scale given by $\Delta := (\log_{10}(9.0) - \log_{10}(3.0)) / 5 \sim 0.095$.
    Thus $\C = \cup_{k=1}^{6} \C_k$, with $\C_k := [-9.0 \times 10^{-(6-k) \Delta}, -9.0 \times 10^{-(5-k) \Delta})$ ($k = 1, \dots, 5$) and $\C_6 := (-\infty, -9.0)$.

    %\item \textbf{Estimation of the occupancy probability of each sub-interval $\C_k$:} The probabilities $p_k := p(\C_k) := \lim_{t \to \infty} p(Z_t \in \C_k)$ are estimated following \cite{mastropietro2025questa}, where the key input of the estimator is the proportion of particles in each set $C_k$ at each time step $t$. We note that the time step $t$  

    \item \label{itm:ours_step_data} \textbf{Rare event generation:} Given the $6$ pairs $(\C_k, \hat{p}_k)_{k = 1, \dots, 6}$, where $\hat{p}_k$ is an estimator of the stationary probability $p(\C_k) := \lim_{h \to \infty} p(Z_h \in \C_k)$, given by e.g. the FV estimator described in Section~\ref{sec:methods_fv}, every realization of a rare wind power change is computed as follows:
    \begin{enumerate}[label=(\roman*)]
        \item an interval $\C_l$ is chosen among $\{\C_k\}_{k = 1, \dots, 6}$ with probability $\frac{\hat{p}_k}{\sum_{j=1}^{6}\hat{p}_j}$;
        
        \item the rare change realization is defined as $y - r$, where $y$ is the upper bound of the selected $\C_l$ interval and $r \sim \mathcal{U}(0, \delta_l)$ if $1 \leq l \leq 5$ (where $\mathcal{U}(a, b)$ is the uniform distribution in interval $(a, b)$, and $\delta_l$ is the length of sub-interval $\C_l$), or $r \sim \mathcal{E}(1.0)$ if $l = 6$ (where $\mathcal{E}$ is the exponential distribution having the enclose rate parameter).
    \end{enumerate}
\end{enumerate}

\subsection{Solution to the multi-stage optimization problem}
The optimization problem to solve using the scenario approach described in Section~\ref{sec:methods} contains $4V = 673\,684$ binary variables where $V = 168\,421$, the number of nodes in the scenario tree. These binary variables define the decisions to be taken for the coal plant operation at each of the $5$ stages of each generated scenario.
The constraints defined in problem \eqref{equ:optim_problem} are thus represented by a sparse matrix with $637\,684$ columns, whose non-zero entries are either $+1$ or $-1$, except for the entries responsible for defining the hard constraint on the coal plant power Eq. \eqref{equ:chance_constraint_as_hard_constraint}, derived from the chance constraint \eqref{equ:optim_problem_constraint_power}.
%converted to a hard constraint as given in associated to the power provided by the coal plant, after conversion into a hard constraint as described in Section~\ref{sec:methods}.

The problem is solved using the standard \texttt{milp} module of the \texttt{scipy.optimize} package in Python. 

\subsection{Evaluation of the solution}
\label{sec:evaluation}
The solution of the optimization problem is evaluated for both \textit{robustness to rare events} and \textit{cost}, on each of the two scenario generation contexts, ``Benchmark'' and ``Biased'', described in Section~\ref{sec:scenario_generation}.

The evaluation is carried out on
%two sets of $100$ replications of
the wind power supply generated for $5$ stages (hours) following a mixture model, as described below.
%As above, the generated wind power is obtained directly from the generated wind speed $Y_h$ by expressing the wind speed in GW.
Note that we distinguish the realized wind power $Y_h$ from the wind power $W_h$ generated as branches of the scenario tree. The value of $Y_h$ will be compared with $W_h$ at each stage in order to find the closest scenario to the realized wind power series.

\begin{figure}
\centering
    \includegraphics[scale=0.55]{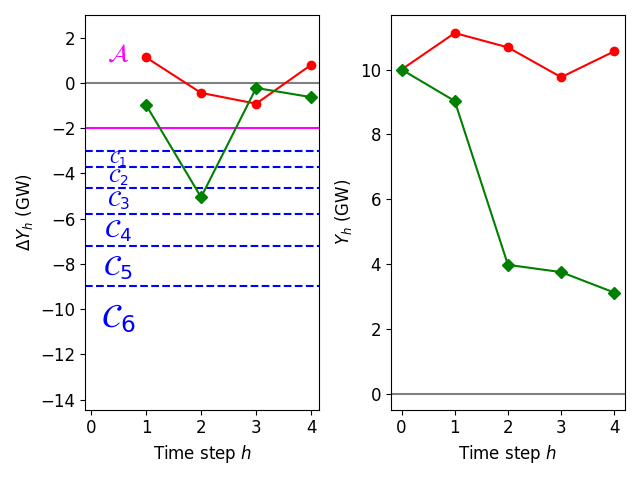}
    \caption{\textbf{Left:} Illustration of two possible realizations of the wind power change process $\Delta Y_h := Y_h - Y_{h-1}, h = 1, \dots, 4$ defined in Eq.~\eqref{equ:wind_power_generation_model}.
    %and of a typical excursion of an FV particle following the dynamics of the $Z_t$ process until it is absorbed at time step $t = 8$.
    The two wind power change realizations are (i) a standard realization with no jumps (red circle symbols), and (ii) a standard realization with one negative jump at time step $h = 2$ with $\Delta Y_2$ falling in the $C_3$ interval (green diamond symbols).
    %The FV particle excursion is shown in cyan square symbols with the absorption event at $t=8$ marked with a cross.
    As an indication of rareness, the sets (intervals) used in the FV simulation responsible for estimating the probability of rare changes are overlaid, namely: the absorption set $\A := [-2.0, +\infty)$ in magenta, and intervals $(\C_k)_{k=1, \dots, 6}$ in blue, whose occupancy probability by the underlying wind power change process $Z_t$ is estimated using the FV estimator and whose disjoint union makes up $\C$ in expression~\eqref{equ:proba_stationary}.
    \\
    \textbf{Right:} Wind power process $Y_h$ for the two realizations of $\Delta Y_h$ on the 5-stage horizon shown on the left, computed as $Y_h = Y_{h-1} + \Delta Y_h, h = 1, \dots, 4$ with $Y_0 = 10$ GW. \\
    %The two wind power change realizations are shown in each plot: (i) in red and circle symbols, a normal realization with no jumps; (ii) in green and diamond symbols, a normal realization with one negative jump at time step $h = 2$ with $\Delta Y_2$ falling in the $C_3$ interval; (iii) in cyan and square symbols, a typical excursion of an FV particle which follows the dynamics of the change process $Z_t$, starting outside but near $\A$ and absorbed at time step $t = 3$ when entering $\A$ (indicated by a cross).
    }
\label{fig:sets_and_realizations}
\end{figure}

The realized $Y_h$ wind power series is obtained using the following mixture model:
\begin{equation}
\label{equ:wind_power_generation_model}
    Y_h = Y_{h-1} + Z_h \ind_{J_h = 0} + R_h \ind_{J_h = 1}, \hspace{0.5cm}\text{for }h = 1, \dots, 4,
\end{equation}
where $Y_0 = 10$ GW (the same initial value for the scenario wind power, $W_0$), $J_h$ is a Bernoulli process with constant probability $q \in [0, 1]$ of $J_h = 1$, drawn independently at each stage $h$, and $R_h$ is a ``rare'' realization of the $Z_h$ process drawn following the procedure described in item~\ref{itm:ours_step_data} of Section~\ref{sec:scenario_generation}.
Two example realizations of $\Delta Y_h$ and $Y_h$ are shown in Figure~\ref{fig:sets_and_realizations}, which highlight the effect on wind power $Y_h$ (right plot) of observing a (rare) negative jump in wind power change $\Delta Y_h$ (left plot), compared to not observing any jump.

The solution's \textit{robustness} is evaluated via demand satisfaction of the solutions across the ``rare'' realizations among the $100$ realizations generated using three different Bernoulli probabilities above, $q = 0\%, 5\%, 10\%$. The solution's \textit{cost} is evaluated by
%setting $q = 0$ (i.e. considering all ``normal'' realizations) and computed as
computing the empirical average cost across all solutions in each case.
%to the $100$ realizations.
A realization is considered ``rare'' if $J_h = 1$ for at least one stage $h$, hence the case $q = 0\%$ works as a baseline having no ``rare'' realizations.

In order to compute these metrics for each of the two evaluation contexts, we proceed as follows:
\begin{enumerate}
    \item \textbf{Obtain the coal plant states for the realization:}
    %In each of the two evaluation contexts,
    For each realization of the wind power supply, $(Y_h)_{h = 0, \dots, 4}$ (red curves in Figure~\ref{fig:realizations}),
    we obtain the coal plant state at each stage $h$ as the subset of the solution $\mathbf{x}$
    %to the optimization problem
    indexed by
    the node in the evaluated scenario tree that is closest to the realized $Y_h$ process in terms of its wind power value, $W_h$ (blue points in Figure~\ref{fig:realizations}).
    %at the given stage $h$.
    The closest node IDs, $n_h\,(h = 0, \dots, 4)$, are found separately for each stage as those minimizing the distance from $W_h$ to the realized $Y_h$ value among the \textit{descendants} of the closest node at the previous stage\footnote{Formally, given $n_0$ the root node ID, we recursively define $n_h := \min_{n \in \mathcal{D}(n_{h-1})} |Y_h - W_h(n)|$ for $h = 1, \dots, 4$, where $W_h(n)$ is the wind power supply at stage $h$ in node $n$ of the scenario tree and $\mathcal{D}(n_{h-1})$ is the set of descendants of node $n_{h-1}$.}.
    %at the first stage the closest node is the root node $n_0$, at the second stage the closest node $n_1$ is found as the one minimizing $|Y_1 - W_1|$, at the third stage the closest node is found as the one ...

    \item \textbf{Compute power and cost of coal plant:} Given the sequence of closest node IDs, $(n_h)_{h = 0, \dots, 4}$, we compute the power provided by the coal plant at each stage $h$ as a function of its state at node $n_h$, $\mathbf{x}(n_h)$, and the \textit{realized} wind power $Y_h$. This power is nonzero only when the plant is operating (i.e. when $x^3(n_h) = 1$), in which case it is equal to $p_h := \min\{p_{\max}, [D_h - Y_h]_+ \}$.
    %, where $x^3(n_h)$ is the value of the binary decision variable $x^3$ at node $n_h$ defining whether the coal plant is operating at stage $h$ (in which case $x^3(n_h) = 1$). The cost associated to the coal plant operation is computed using \eqref{equ:cost}.
\end{enumerate}
%\textcolor{red}{This is a fair amount of detail for a reader to digest, perhaps you can make a tikz illustration?}

\subsection{Results}
Table~\ref{tab:results} summarizes the cost and robustness analyses described in the preivous section for the two scenario generation strategies (Benchmark and Biased) and the three Bernoulli probability of a rare negative change in wind power at each stage.
The realizations and closest scenarios behind those numbers are shown in Figure~\ref{fig:realizations}.
%shows all the realizations behind the results presented in Table~\ref{tab:results}, grouped by the two scenario generation contexts (``BM'' and ``Biased'') and
%to evaluate the solution's \textit{robustness} and \textit{cost} for the two scenario generation contexts considered.
%the three different Bernoulli probabilities $q$ of rare changes in wind power.

From the table, we observe the following:
\begin{enumerate}
    \item The average costs increase as the proportion of rare instances increases.
    \item The observed average cost of coal plant operation across realizations is always larger in the Biased scenario tree approach than in the Benchmark (``BM'') approach, although not larger than twice as much. As counterpart, the Biased approach shows a $100\%$ satisfaction of demand, which is not the case in the Benchmark approach, unless no rare instances happen ($q = 0\%$ column on the left).
    \item Naturally, the demand unsatisfaction in the Benchmark approach increases as the proportion of rare instances increases. We note also that most of the unsatisfied power demand comes from the rare instances (as deduced from the similar percentages in items $3$ and $5$ of the table). 
\end{enumerate}  

From the figure, and as expected, we observe that the
bulk
%density
of blue curves (representing the closest scenario to each realization in red) is
%higher in the region of normal wind power supply in the ``BM'' context (left column) than in the ``Biased'' context (right column).
biased towards lower wind power regions in the Biased approach compared to the Benchmark.
%In the latter, on the contrary, the density of blue lines is distributed more equally between the ``normal'' region and the ``rare'' region represented by low wind power values. In fact, 
In fact, thanks to the negative bias in the Biased scenario approach,
%more than the benchmark case, thus providing
the optimization model is able to better learn
%with a higher chance of knowing
what the best decision for the coal plant is under rarely observed shortfall instances.

In sum, even though the Benchmark approach yields about half the average cost for plant operation compared to our Biased approach, it is not as robust as the latter which is able to satisfy demand in all instances.
As the resolution of the approximated problem gets closer to the actual infinite dimensional optimization by increasing the number of scenarios (albeit at the cost of more computational power), the gap between the respective average costs should decrease. However, the higher risk of unsatisfied demand in the Benchmark approach will be difficult to overcome considering the existing considerable leap with the Biased approach.
%given it is based on an event occurrence rather than on a continuous quantity.

%\textcolor{red}{Some commentary on the results: \textit{the cost shortfall is (only) double}, suggesting more computational power can be enough to mitigate it, while there is a substantial leap in the risk of unsatisfied demand for rare instances in the baseline}
%The metrics related to the robustness evaluation in each context are reported in Table~\ref{tab:results}.

\section{Conclusion}\label{sec:conc}
We have presented a methodology to take into account rare occurrences of a process to obtain optimal solutions of a stochastic optimization problem. We illustrated the approach in the context of electricity grids powered by renewable energy sources (such as wind), with alternative backup sources (such as coal) to compensate for shortfalls due to the high volatility in renewable energy generation. 

The numerical results authoritatively demonstrate that the use of scenario generation together with rare event sampling bias using Fleming-Viot can eliminate energy supply short fall risk, with relatively moderate tradeoffs as far as average cost. This validates the potential utility of the Fleming-Viot approach for rare event sampling in the context of robust stochastic programming. We expect that more comprehensive simulations that include a greater array of choices for energy generation sources, and over longer horizons, would yield similar results. One potentially interesting approach would be to have two layers of robustness for rare and very rare scenarios by considering both  classical coal plants (for the very rare scenarios) and the cleaner combined cycle gas turbine (CCGT) plants (for the rare scenarios). 

Future work can include modeling Demand as a stochastic process and more physically faithful wind speed and wind power modeling. More formally, an exploration of scaling with scenario tree size would be interesting and informative to operators. Specifically, the two approaches become theoretically identical as the number of scenarios approaches infinity, but we conjecture that the rate of convergence in average cost and rare event robustness for the two methods is distinct.

\renewcommand{\arraystretch}{1.2}
\setlength{\tabcolsep}{4pt}
\begin{table}%[!ht]
    \centering
    %\tiny
    \scriptsize
    %\footnotesize
    %\normalsize
    \begin{tabular}{l!{\vrule width 2pt}c|c!{\vrule width 1.2pt}c|c!{\vrule width 1.0pt}c|c}
        \begin{tabular}{@{}l@{}}\textbf{Performance metric} \\
        on $n = 100$ realizations \\
        (lower is better)
        \end{tabular} & \multicolumn{2}{c!{\vrule width 1.0pt}}{\begin{tabular}{@{}c@{}}$q = 0\%$ \\ ($n_r = 0$ rare)\end{tabular}} & \multicolumn{2}{c!{\vrule width 1.0pt}}{\begin{tabular}{@{}c@{}}$q = 5\%$ \\ ($n_r = 19$ rare)\end{tabular}} & \multicolumn{2}{c}{\begin{tabular}{@{}c@{}}$q = 10\%$ \\ ($n_r = 36$ rare)\end{tabular}} \\
        %\hline
        \cline{1-7}
        %& $q = 0\%$ & 0 ``rare''
        %%& 3 ``rare'' & 
        %& $q = 5\%$ & 19 ``rare''
        %& $q = 10\%$ & 36 ``rare'' \\
        \multicolumn{1}{r!{\vrule width 2pt}}{Scenario tree method $\rightarrow$} & \textbf{BM} & \textbf{Biased}
        %& \textbf{BM} & \textbf{Ours} 
        & \textbf{BM} & \textbf{Biased} & \textbf{BM} & \textbf{Biased} \\
        
        \noalign{\hrule height 2pt}
        \begin{tabular}{@{}l@{}}(1) Avg. observed cost\end{tabular} & \textbf{12.9} & 22.7 
        %& 12.9 & 23.2
        & \textbf{19.9} & 32.2
        & \textbf{26.5} & 42.9 \\ %17.2 & 22.3 \\
        
        \noalign{\hrule height 1pt}
        \begin{tabular}{@{}l@{}}(2) \% Realizations with \\ some unsatisfied demand\end{tabular} & 0.0\% & 0.0\%
        %& 2.0\% & \textbf{0.0\%} 
        & 6.0\% & \textbf{0.0\%} & 13.0\% & \textbf{0.0\%} \\ %16.0\% & 1.0\% \\
        \hline
        \begin{tabular}{@{}l@{}}(3) \% Unsatisfied power \\ demand in all realizations\end{tabular} & 0.0\% & 0.0\%
        %& 30.9\% & \textbf{0.0\%} 
        & 23.7\% & \textbf{0.0\%} & 30.4\% & \textbf{0.0\%} \\ %26.7\% & 0.3\% \\

        \noalign{\hrule height 1pt}
        \begin{tabular}{@{}l@{}}(4) \% \textbf{Rare} realizations \\ with some unsatisfied \\ demand\end{tabular} & --- & --- 
        %& 66.7\% & \textbf{0.0\%} 
        & 31.6\% & \textbf{0.0\%} & 36.1\% & \textbf{0.0\%} \\ %44.4\% & 2.8\% \\
        \hline
        \begin{tabular}{@{}l@{}}(5) \% Unsatisfied power \\ demand in \textbf{rare} \\ realizations\end{tabular} & --- & --- 
        %& 97.5\% & \textbf{0.0\%} 
        & 24.6\% & \textbf{0.0\%} & 30.6\% & \textbf{0.0\%} \\ %26.9\% & 0.3\%
    \end{tabular}
    \caption{Evaluation metrics (lower is better) measuring the ability of a benchmark strategy (``BM'') and our proposed strategy (``Biased'') for controlling the operation of a backup coal plant in mitigating the effect on demand satisfaction of power supply shortfall from renewable sources.
    The metrics are based on $100$ realizations of the wind power supply over a 5-stage period using Eq.~\eqref{equ:wind_power_generation_model} for three different values of the Bernoulli probability $q$ (across columns), driving negative jumps in wind power (i.e. rare negative changes $Y_h - Y_{h-1}$). A realization is considered ``rare'' if at least one consecutive change in wind power is the result of a jump.
    %in Section~\ref{sec:results}.
    %out of which 36 realizations have at least one change in wind power supply that is considered rare, thus classifying the whole realization as ``rare''.
    }
\label{tab:results}
\end{table}

\begin{comment}\label{Comment this out if we show the results on different q values in a single table above}
\begin{table}[!ht]
    \centering
    \begin{tabular}{l|c|c}
        \begin{tabular}{@{}l@{}}\textbf{Performance metric (lower is better)} \\ \footnotesize (on $n=100$ ``normal'' realizations)\end{tabular} & \textbf{Benchmark} & \textbf{Ours}  \\
        \hline
        \hline
        Expected cost of optimal decision (solution) & 0.66  & 15.6 \\
        \hline
        Average observed cost of realizations & 12.9 & 22.7 \\
        \% Realizations with some unsatisfied demand & 0.0\% & 0.0\% \\
        \% Unsatisfied power demand in all realizations & 0.0\% & 0.0\% \\
    \end{tabular}
    \caption{Evaluation metrics measuring the ability of a benchmark strategy (``BM'') and our strategy (``Ours'') for controlling the operation of a backup coal plant in mitigating the effect on demand satisfaction of \textbf{``normal''} power supply shortfall by renewable sources. The metrics are based on $100$ realizations of the wind power supply over a 5-stage period, as described in Section~\ref{sec:results}.}
\label{tab:results_cost}
\end{table}
\end{comment}

% Width factor for each plot
\newcommand\factorsub{0.24}
\newcommand\factor{1.0}
\begin{figure}%[!ht]
    \centering
    %\centering Benchmark \qquad \qquad \qquad \qquad Biased \\
    \footnotesize Probability of rare wind power change, $q = 0\%$ \\
    \begin{subfigure}[b]{\factorsub\textwidth}
        \includegraphics[width=\factor\linewidth]{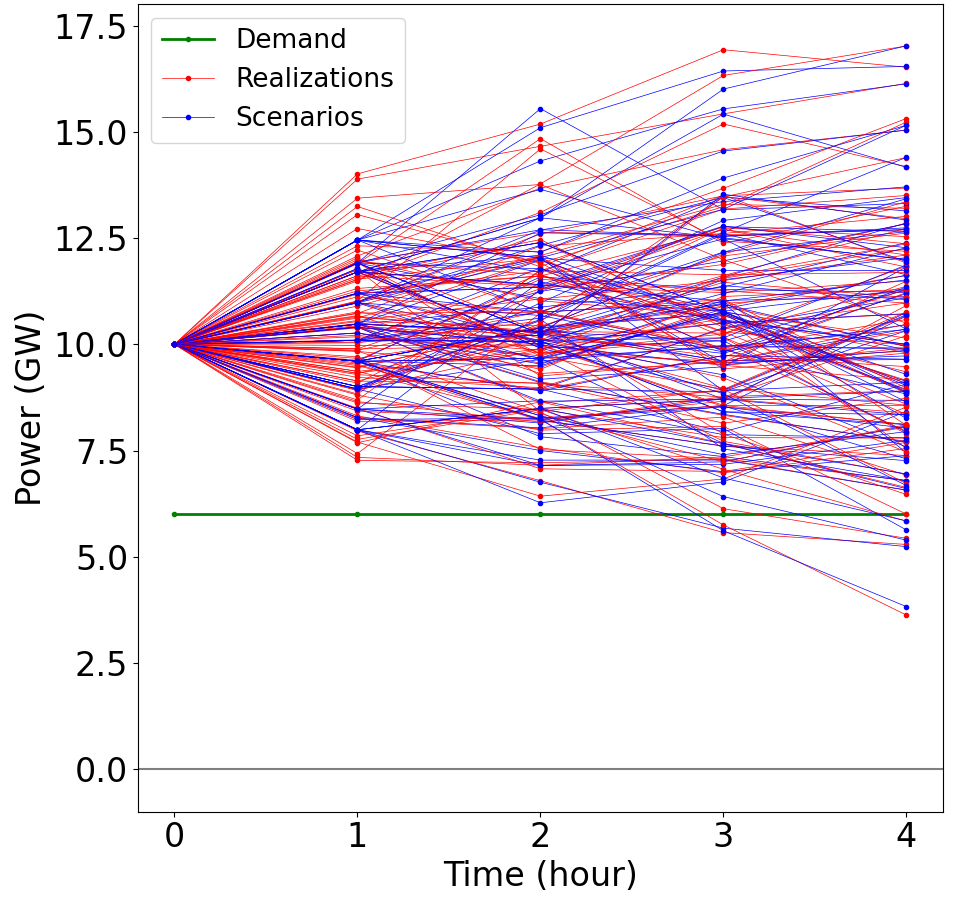}%{figures/MPC_Evaluation_ARPositive_100Realizations10p_OnNORMALScenarios.png}
        %\caption{$q = 0\%$, Benchmark scenario tree}
    \end{subfigure}
    \begin{subfigure}[b]{\factorsub\textwidth}
        \includegraphics[width=\factor\linewidth]{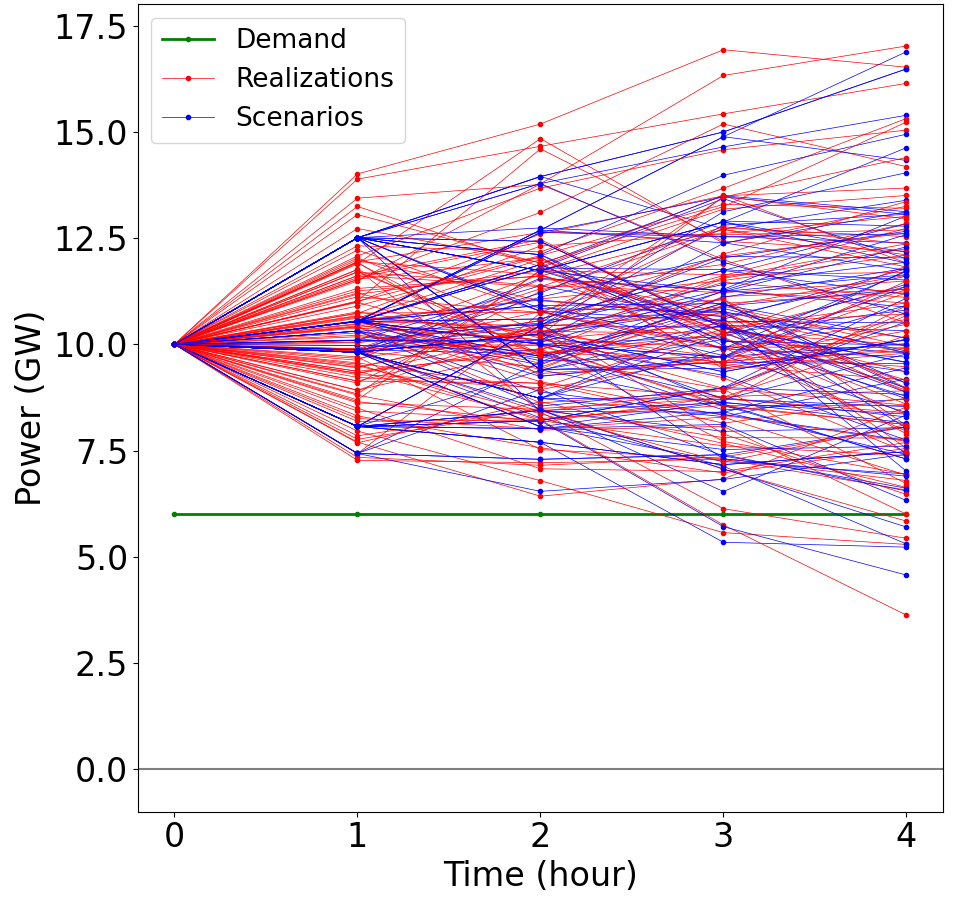}%{figures/MPC_Evaluation_ARPositive_100Realizations10p_OnBIASEDScenarios.png}
        %\caption{$q = 0\%$, Biased scenario tree}
    \end{subfigure}
    \footnotesize Probability of rare wind power change, $q = 5\%$ \\
    \begin{subfigure}[b]{\factorsub\textwidth}
        \includegraphics[width=\factor\linewidth]{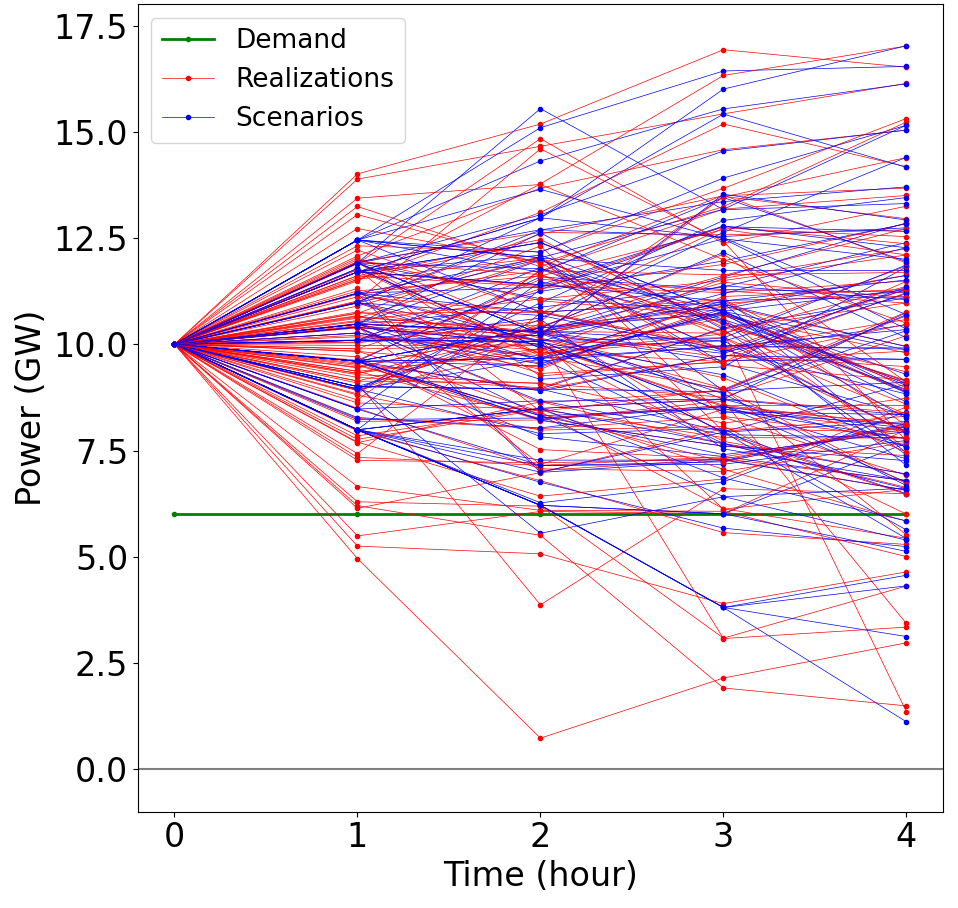}%{figures/MPC_Evaluation_ARPositive_100Realizations10p_OnNORMALScenarios.png}
        %\caption{$q = 5\%$, Benchmark scenario tree}
    \end{subfigure}
    \begin{subfigure}[b]{\factorsub\textwidth}
        \includegraphics[width=\factor\linewidth]{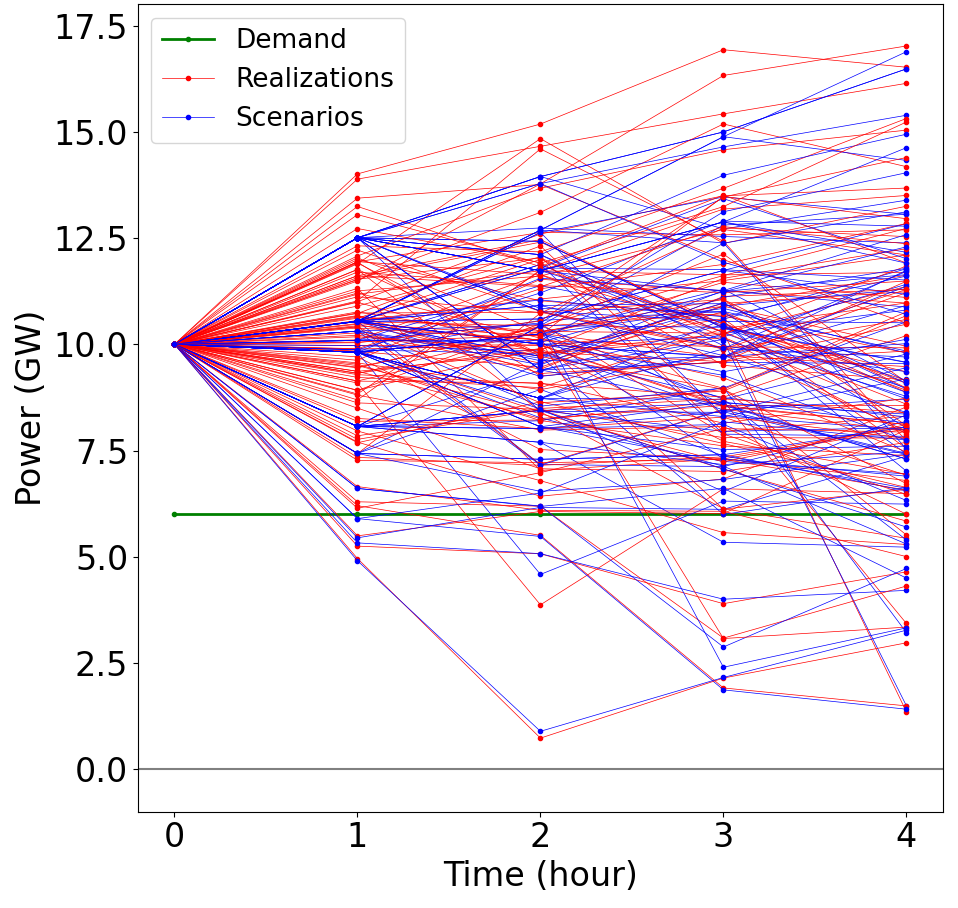}%{figures/MPC_Evaluation_ARPositive_100Realizations10p_OnBIASEDScenarios.png}
        %\caption{$q = 5\%$, Biased scenario tree}
    \end{subfigure}
    \footnotesize Probability of rare wind power change, $q = 10\%$ \\
    \begin{subfigure}[b]{\factorsub\textwidth}
        \includegraphics[width=\factor\linewidth]{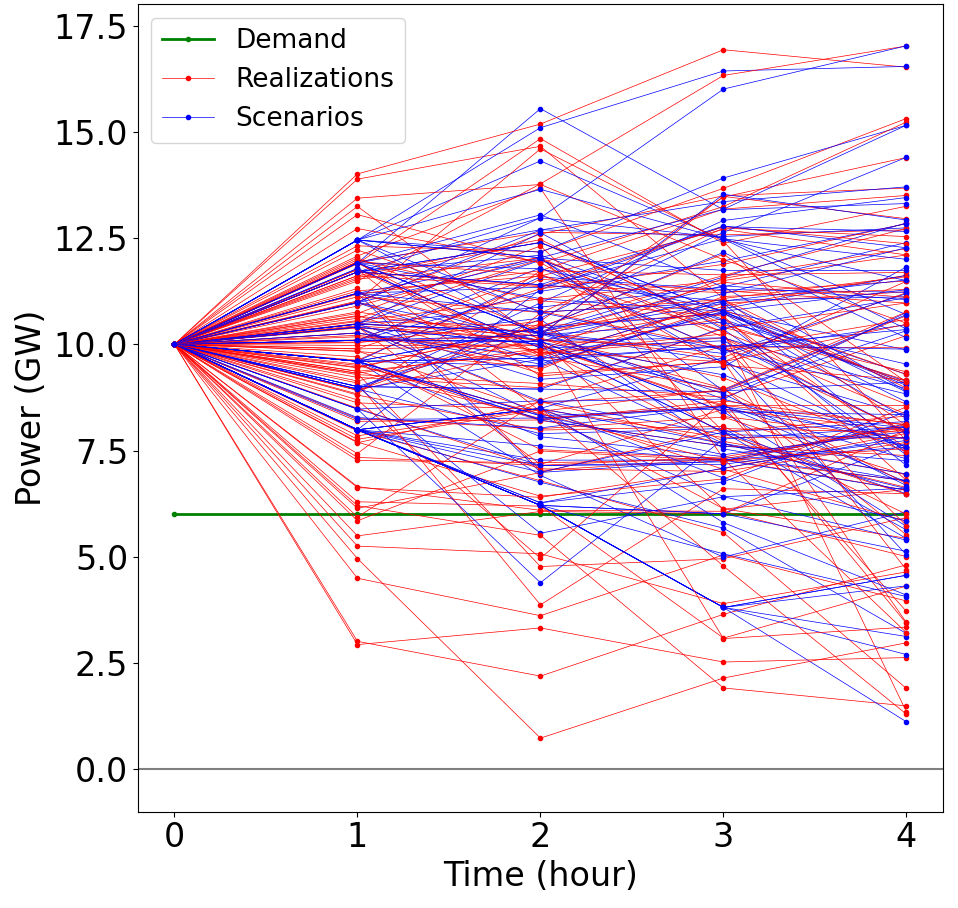}%{figures/MPC_Evaluation_ARPositive_100Realizations10p_OnNORMALScenarios.png}
        %\caption{$q = 10\%$, Benchmark scenario tree}
    \end{subfigure}
    \begin{subfigure}[b]{\factorsub\textwidth}
        \includegraphics[width=\factor\linewidth]{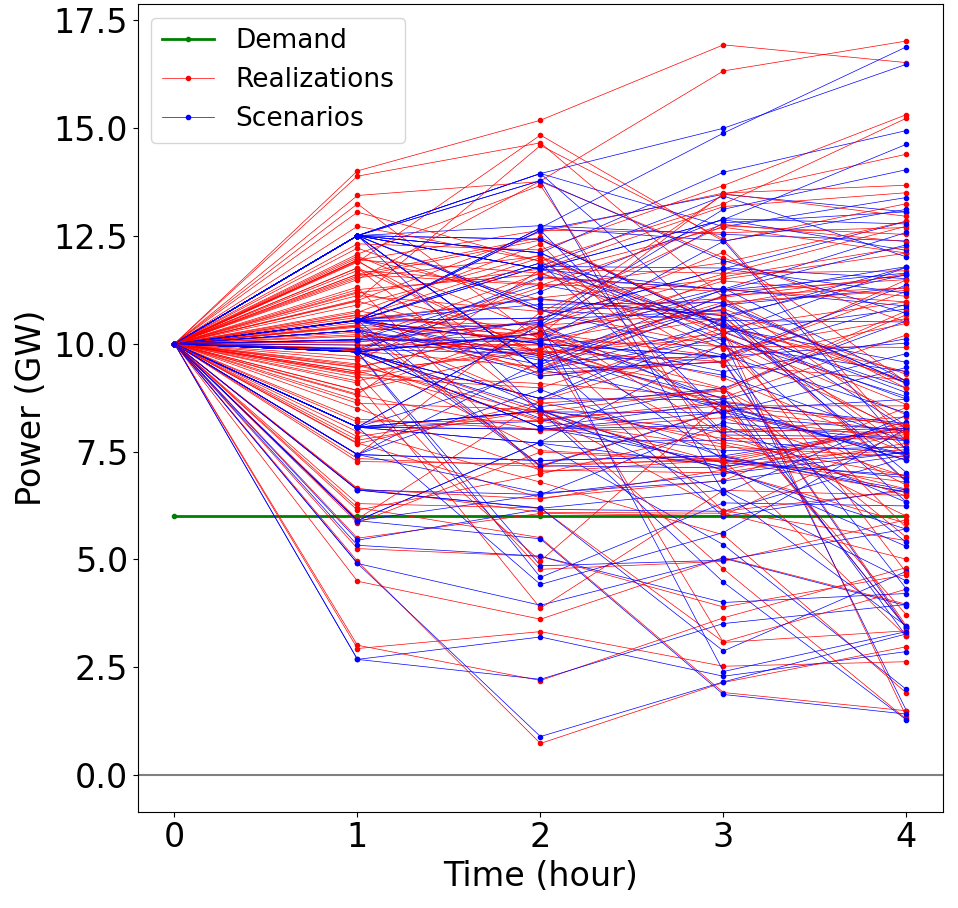}%{figures/MPC_Evaluation_ARPositive_100Realizations10p_OnBIASEDScenarios.png}
        %\caption{$q = 10\%$, Biased scenario tree}
    \end{subfigure}    
    \caption{Wind power realizations (red) and their closest scenario (blue) in each of the two scenario generation strategies considered in the resolution of the 5-stage optimization problem \eqref{equ:optim_problem}: Benchmark (left column) and Biased (right column) (see Section~\ref{sec:scenario_generation}). Each row corresponds to a different Bernoulli probability $q$ of generating a rare wind power change at each stage using Eq.~\eqref{equ:wind_power_generation_model} in Section~\ref{sec:evaluation}.}
\label{fig:realizations}
\end{figure}

\twocolumn

\bibliographystyle{plain}

\end{document}